\documentclass{article}
\usepackage{dcolumn}
\usepackage{multirow}
\usepackage{bm}

\usepackage[pdftex]{graphicx}
\usepackage{amsmath}
\usepackage{amsfonts}
\usepackage{lipsum}
\usepackage{multirow}
\usepackage{authblk}
\usepackage{fancyhdr}
\usepackage{enumitem}
\usepackage{amssymb}
\usepackage{color}
\usepackage{titletoc}
\usepackage{cite}
\usepackage{titletoc}

\usepackage[utf8]{inputenc}
\usepackage[T1]{fontenc}
\usepackage{mathptmx}

\usepackage{hyperref}
\usepackage{url}

\begin{document}

\title{\bf Reconstructing dynamical networks via feature ranking}

\author[1,2]{Marc G.\ Leguia}
\author[1,3]{Zoran Levnaji\'c}
\author[4,3]{Ljup\v co Todorovski}
\author[3]{Bernard \v Zenko}

\affil[1]{\small Faculty of Information Studies in Novo Mesto, Slovenia}
\affil[2]{\small Department of Communication and Information Technologies, Universitat Pompeu Fabra, Barcelona, Spain}
\affil[3]{\small Department of Knowledge Technologies, Jo\v zef Stefan Institute, Ljubljana, Slovenia}
\affil[4]{\small Faculty of Public Administration, University of Ljubljana, Slovenia}

\maketitle


\begin{abstract}
        Empirical data on real complex systems are becoming increasingly available. Parallel to this is the need for new methods of reconstructing (inferring) the structure of networks from time-resolved observations of their node-dynamics. The methods based on physical insights often rely on strong assumptions about the properties and dynamics of the scrutinized network. Here, we use the insights from machine learning to design a new method of network reconstruction that essentially makes no such assumptions. Specifically, we interpret the available trajectories (data) as \textit{features}, and use two independent feature ranking approaches -- Random Forest and RReliefF -- to rank the importance of each node for predicting the value of each other node, which yields the reconstructed adjacency matrix. We show that our method is fairly robust to coupling strength, system size, trajectory length and noise. We also find that the reconstruction quality strongly depends on the dynamical regime.
\end{abstract}


\maketitle


\section{Introduction}

A foremost problem in modern network science is how to reconstruct (infer) the unknown network structure from the available data~\cite{newman,costa,estrada,easley,alb,mason}. Namely, while the functioning of real complex networks can often be to some degree observed and measured, their precise structure (organization of connections among the nodes) is almost never accessible~\cite{liu}. Yet understanding the architectures of real complex networks is key, not just for applied purposes, but also for better grasping their actual functioning~\cite{hebert,ruths}. For this reason, the topic of developing new and efficient methods of network reconstruction gained ground within network science~\cite{timme1}.

This problem is in literature formulated in several ways. Typically, one considers the nodes to be individual dynamical systems, with their local dynamics governed by some difference or differential equation~\cite{mason}. The interaction among these individual systems (nodes) is then articulated via a mathematical function that captures the nature of interactions between the pairs of connected nodes (either by directed or non-directed links). In this setting, the problem of network reconstruction reduces to estimating the presence/absence of links between the pairs of nodes from time-resolved measurements of their dynamics (time series), which are assumed available. It is within this formulation that we approach the topic in this paper, i.e., we consider the structure of the studied network to be hidden in a "black box", and seek to reconstruct it from time series of node dynamics (i.e., discrete trajectories). 

Within the realm of physics literature, many methods have been proposed relying on above formulation of the problem, and are usually anchored in empirical physical insights about network collective behavior~\cite{timme1,timme2,timme3}. This primarily includes synchronization~\cite{alex}, both theoretically~\cite{albert,luce} and experimentally~\cite{blaha,kralemann}, and in the presence of noise~\cite{stankovski1}. Other methods use techniques such as compressive sensing~\cite{sensing,grebogi} or elaborate statistics of derivative-variable correlations~\cite{us1,marc}. Some methods are designed for specific domain problems, such as networks of neurons~\cite{rok,arkady} or even social networks~\cite{kocarev}. There are also approaches specifically intended for high-dimensional dynamical system, mostly realized as phase space reconstruction methods \cite{PhysRevLett.45.712,PhysRevA.45.3403,Ma2018,doi:10.1063/1.5023860}. While many methods in general refer to non-directed networks, some aim specifically at discerning the direction of interactions (infer the 'causality network'). One such method is termed Partial Mutual Information from Mixed Embedding (PMIME~\cite{Koutlis2016a}) and will be of use later in this work.

However, a severe drawback of the existing physical reconstruction paradigms is that the empirical insights that they are based on are often translated into very strong assumptions about our knowledge of the system. For example, these methods often require the knowledge of not just the mathematical form of the dynamical model, but also the precise knowledge of the interaction function(s)~\cite{us1,marc}. Similarly, some methods require the possibility to influence the system under study, for example by resetting its dynamics or influencing it in other ways~\cite{us2,timme2}. Other methods make assumptions about the dynamical nature of the available trajectories (data), e.g., their linearity. Another family of methods require that the mathematical form of interaction function is sparse in non-zero terms~\cite{sparsereg}. While such data-driven methods are elegant and in principle efficient~\cite{Schelter2006,Jachan2009,Kugiumtzis2013,Rubido2014a,tirabassi2015inferring,bianco2016successful,Koutlis2016a,Leguia2019}, these methods often require long data-sets and/or the implicit assumptions about the signals that can be limiting to their usage in some situations of practical interest. In fact, latest results emphasize the importance of model-free reconstruction methods~\cite{sparsereg,timme2,timme17nc}, which is the context of our present contribution.

Moreover, \emph{relevance network approach (RNA)} follows the statistical perspective of the network reconstruction task and is often used for inferring gene regulatory networks from expression data~\cite{relevance,Hemp11a}. There, the decision on the link presence or absence is based on the pairwise correlation between the time series observed in nodes. There are many different variants of RNA, each corresponding to a different measure of association between time series. The commonly used association measures include Pearson correlation coefficient and entropy-based mutual information. When comparing with the feature-ranking approach presented here, we need to emphasize the fact that our approach takes a multi-variate view on the correlation, since the predictive model for each node is built taking into account the (potential) influence of all the other nodes in the network. In contrast, RNA takes a uni-variate stance on correlation, by measuring it against a time series in each network node separately. Additionally, each measure of pairwise correlation used in RNA often assumes a known influence model, e.g., Pearson correlation coefficient (one of the measures used in RNA) assumes that the interdependence between the given pair of nodes is linear.

On a parallel front, the problem of network reconstruction (inference) has been tackled also by computer scientists, specifically within the field of machine learning~\cite{pat,vapnik,hastie}. The key difference between physics and machine learning paradigms is that the latter make very little or no assumptions about our knowledge of the system, i.e., approach the problem in full generality. And while such methods can be resource-demanding, they are actually more applicable in real scenarios of many domain sciences, where a given network is to be reconstructed with basically no prior knowledge about its internal details~\cite{toussaint,zanin,wen}. Various methods have been developed along these lines, based on techniques of equation discovery~\cite{dzeroski} and symbolic regression~\cite{lipson,sparsereg}. The symbolic regression methods based on sparse regression~\cite{sparsereg} and block-sparse regression~\cite{timme17nc} has been successfully used to reconstruct small-scale networks~\cite{sparsereg-nets}. Note however, that these methods not only aim at reconstructing the network structure, but also infer the mathematical models of the interactions leading to complete reconstruction of the network dynamics. The later render these reconstruction methods computationally expensive. 
Scalability tends to be a problem of network reconstruction methods in general, and existence of computationally more efficient approaches would significantly increase the number of practical problems that could be tackled.

Among the core techniques in machine learning is \textit{supervised learning}: one tries to learn from the observation data how does a dependent variable (\emph{target}) depend on a given set of independent variables (\emph{features}). To this aim, one searches for a (predictive) mathematical model that is to capture this dependence. This model can also be used to predict the value of the target given the values of the features. In such a model, not all features will play the same role---the target variable will in general depend more on some features than on others. We can therefore rank the features according to their influence on the target, and this is what machine learning literature calls \emph{feature ranking}~\cite{Guyon:2003}. There is a range of different feature ranking methods, such as RReliefF~\cite{rrelieff} and Random Forest~\cite{BREIM01}, and with a ranking produced by one of these methods, one can improve the learned model in several ways. The simplest of them is to ignore the features (independent variables) with low ranks, as they have little or no influence on the target. Such features are often complicating the model without contributing to its accuracy. In fact, a simplified model without such features can be even more accurate (due to a phenomenon called over-fitting)~\cite{Guyon:2003}. Crucial then is to set the best threshold on which features to ignore and which to keep in order to obtain the most accurate model.  

In this paper we propose a new method of reconstructing a dynamical network of physical interest from discrete time series of node dynamics. In contrast to the usual formulations of this problem in physics literature, we here build our reconstruction method on the concept of feature ranking. Specifically, we treat the dynamical state of a given node as the target variable, while the previous states of all other nodes are treated as features. We use the dynamical data to quantify how much each feature influences the target, and compute the feature ranking accordingly. Some features will have a strong influence on the target, so it is reasonable to assume that the corresponding nodes are linked to the studied node. It is also safe to assume that low ranked features (nodes) are not connected to the studied node: Of course, when using this method in practice, one has to carefully select the value of the threshold by observing the sensitivity-specificity trade-off. In the evaluation of our method, we study the reconstruction performance for all candidate threshold values using Receiver Operator Characteristics (ROC) curve.

Note that in the formulation of our method we made \emph{no assumptions} about the knowledge of the interaction functions or the dynamical equations of network dynamics. Therefore, our method relies solely on the time series and their properties such as length and possible presence of observational noise (we assume that time series are coming from some empirical measurement/observation of the system). 

The rest of the paper is organized as follows. In the next Section we first explain some basic concepts from machine learning and feature ranking, and then explain and motivate our method. In Section Results we illustrate the performance of our method using several examples and study its response to diverse properties of the system. We close the paper with the discussion of our findings and limitations of our method, emphasizing the potentials for practical use.


\section{The reconstruction method}

In this section we explain our reconstruction method. For clarity we build it step by step, first explaining the relevant concepts of its machine learning background.


\paragraph{Machine learning, features and feature ranking.} Machine learning studies algorithms whose performance improves with 'experience'~\cite{Mitchell1997}. Such improvement is typically gained by making the algorithm 'learn' from that experience, which comes in form of many examples of data~\cite{HTF06,WITT15}. To 'learn' means to look for patterns in the data and extract them: For example, by making a Fourier decomposition of various sound signals, one can ``learn" to differentiate between human speech and birdsong. Machine learning can be seen as an approach to data-driven modeling suited for circumstances when our knowledge about the studied system is limited. This is the core reason why machine learning is being increasingly used in a variety of scientific disciplines, ranging from from medicine and biology~\cite{furey2000support,shipp2002diffuse,guyon2002gene}, to stock market analysis~\cite{huang2005forecasting}, text classification~\cite{TSKD01,sebastiani2002machine} and image identification~\cite{rosten2006machine}.

Physics community has over the past decade recognized this ability of machine learning, which triggered an array of novel results in diverse fields of physics~\cite{schoenholz,seko,hentschel,yi,toussaint}, including complex networks~\cite{wen} and dynamical systems~\cite{zanin}. In particular, machine learning was also used to formulate the network reconstruction problem for several domain sciences~\cite{huynh-thu2010plosone,decoupling}.


In the most common setting of supervised learning, an algorithm uses existing examples of data as inputs, and produces a set of patterns or a \textit{predictive model} as the output. Examples of data are typically given in the \emph{attribute-value representation}~\cite{ziarko1996attval}, 
which means that each data example is described via series of values of \emph{attributes} (in machine learning also called \emph{features} or \emph{independent variables}). Hence, one can use the input data to create a predictive model describing how the target variable depends on the features. The model can be then used to predict the value of the target, given any values of the features, even ones not included in the training data. Furthermore, the model can be used to determine the importance of features or feature ranks.

To illustrate the idea of feature ranking, say we are given the equation
\begin{equation}
  y = f(x_{1},x_{2},...) = x_{1}^2 + x_{2} + 2,
  \label{ex:1}
\end{equation}
and let us assume that function $f$ is not known, but it is known that $y$ depends on several variables $x_i$. In other words, $y$ is the target variable and $x_i$ are the features. This type of a task is referred in machine learning as a regression task and it is being solved using a machine learning algorithm $M$, i.e., $f \approx \hat{f} = M(D)$, where $D$ is the data set and $\hat{f}$ is the prediction model that for any given observation $(x_1,x_2,x_3)$ can be used to predict the value of $y$, $\hat{y} = \hat{f}(x_1,x_2,x_3)$. Suppose now that we are given the following list of values (measured with some observation error). In other words, we are given a data set $D$ consisting of $L$ attribute-value tuples $(x_1, x_2, x_3; y)$:
\medskip\noindent
\begin{center} \begin{tabular}{@{}lrrrr@{}} 
  example & $x_1$ & $x_2$ & $x_3$ & $y$\\\hline
  \#1 & 2.1 & 1.9 & 2.3 & 8.0 \\
  \#2 & 4.7 & 0.7 & 5.3 & 27.4 \\
  \ldots \\
  \#$L$ & 10.6 & 7.9 & 4.5 & 114.8\\
\end{tabular} \end{center}
\medskip
\noindent and we want to reconstruct (or infer) 
$f$. Note that in the data we also have the feature $x_3$, which actually does not influence $y$, but we assume not to know that a-priori. This situation is very common in various scientific domains, as the inspected system is often poorly understood, and the only available data is collected via features that \emph{may or may not} influence the target variable.


One example of a machine learning algorithm $M$ for regression is \emph{Random Forest}~\cite{BREIM01}. A Random Forest model is an ensemble of piece-wise constant models, regression trees~\cite{BFJRCR84}, where model segments correspond to intervals of feature values: The algorithm for learning regression trees uses training data to approximate the optimal splits of the data space into segments with a constant value of the target (regression tree is a hierarchical splitting of the feature space into segments). Each tree in the Random Forest ensemble is learned on a random sample of the learning data set $D$, and each split in the tree is chosen from a random sample of features $x_i$. The prediction $\hat{y}$ of the ensemble is the average of predictions of all the trees. So, each random sample of input data gives a new tree (an independent splitting scheme), which we average over. While learning a single tree is prone to over-fitting the training data, the ensemble of trees is proven to be more robust, thus leading to accurate predictive models.

The Random Forest machine learning algorithm can be also used for feature ranking. One can compare the prediction error of (i) the Random Forest model learned on the learning data set with the prediction error of (ii) the Random Forest model learned on the randomized training data, where the values of the feature of interest are being randomly permuted between the data points. Intuitively, if the errors of the two models differ a lot, then the importance of the feature is high: Note that in this case, the randomly permuted values of the feature cause the model error increase, hence the feature contributes a lot to the model accuracy. And vice versa: The feature importance is small, if the observed difference is low. 


Another example of an algorithm for regression is \emph{Nearest Neighbor}~\cite{aha1991}. Given a data point $x$, the Nearest Neighbor algorithm finds its nearest neighbors in the learning data set $D$ (with respect to the values of the features) and then predicts the target value of $x$ as an average of the target values of the nearest neighbors with respect to a distance measure (e.g. Euclidean) in the feature space. \emph{RReliefF}~\cite{Kononenko1997} is an extension of the simple nearest neighbor idea for feature ranking. It ranks the importance of features based on the detected differences between nearest neighbor input data example pairs: If there is a feature value difference in a pair with the similar target value, the feature importance is decreased. In contrast, if there is a feature value difference in a pair with dissimilar target values, the feature importance is increased.

Let us assume now that we applied a feature ranking algorithm $R$ (such as Random Forest or RReliefF) on the above data set $D$ and obtained the following feature ranking or ranking scores (values are illustrative):
$$
R(D) = (F_1,F_2,F_3) = (12.3,2.5,0.2),
$$
\noindent where each $ F_i $ denotes the importance of the feature $ x_i $. The exact values of the ranking scores are not important, what is important are their relative values. In this case $x_1$ has the largest score, which means that it is ranked as the most important feature for the target variable $y$. $x_3$, on the other hand, has the lowest score and its influence on the value of $y$ is small, if such influence exists at all (from Eq.~\ref{ex:1} we know that it actually does not). This ranking can now be used as input for modeling $y$ with one of the standard regression methods: Instead of fitting it on all three features, we fit only on $x_1$ and $x_2$. However, in practice, deciding where to draw the line and what features to ignore is far from trivial and often depends on the particularities of the problem at hand. In machine learning the issue of identifying the relevant features is a classic problem in its own right called \emph{feature selection}, and can be studied via several approaches, including feature ranking.

\paragraph{Our reconstruction method.} Armed with above insight, we now proceed to our reconstruction method. As already mentioned above, feature ranking methods can be naturally applied to the problem of reconstructing a dynamical network from the observations (time series measurements) of its node dynamics. Assuming that the state of a selected node of a dynamical network represents a target, and that its state is influenced by the connected network nodes, which represent features, one can define a supervised learning problem for this node: Learn a regression model for predicting the state of the selected node from the states of all the other nodes in the network. Note that our aim here is not the predictive model, instead, we are interested in the feature ranking only: we can actually rank the importance of the other network nodes to the selected one, because a highly ranked node is likely to be connected to the selected node. We now only need to repeat this procedure for all the nodes and we can reconstruct the entire network structure. 

Note that this articulation of the reconstruction problem includes \emph{no assumptions} about the network dynamical model or general properties of the dynamics. 

Let us now formally present our method. Although we have developed it independently, it is very similar to the method presented in~\cite{huynh-thu2010plosone}. We start with a general network with $N$ nodes. The state of a node $i$ at time $t+1$ is $x_i(t+1)$, and it dynamics is influenced by the states of all nodes connected to $i$ at some earlier time $t$:
\begin{equation}
  x_{i}(t+1) = f_{i}(x_{1}(t),x_{2}(t),\dots,x_{N}(t));\quad i=1\dots N  .
  \label{eq:featureranking}
\end{equation}
Note that this is the most general possible formulation of the network dynamics: 
Each node's behavior is influenced by unknown node-specific interaction function $f_i$, dependent on all other nodes. We assume total observability of the system, meaning we have access to the trajectories of all nodes at all times. 

When reconstructing the network, the interaction function $f_i$ is not known, but we can use the observation data to model it. The observation data consists of state trajectories $(x_{i}(1),x_{i}(2), \dots, x_{i}(L))$ for all the network nodes. The Eq.~\ref{eq:featureranking} therefore represents our regression modeling problem for node $i$, where the state variable $x_i(t+1)$ is the target and state variables $x_j(t); j = 1\dots N$ are the features. From the observation data we construct the training data with $L-1$ examples
\begin{equation}
  D_i = \bigcup_{t=1}^{L-1}\; \big(x_1(t),x_2(t),\dots,x_N(t); x_i(t+1)\big) ,
\end{equation}
and a suitable machine learning algorithm for regression $M$ could be used to compute the approximation $\hat{f}_i$
\begin{equation}
    f_i \approx \hat{f}_i = M(D_i).
\end{equation}
%
%
However, we are not really interested in solving these $N$ regression problems, we only perform feature ranking for them, since it is these rankings that contain information on the connections among the nodes. In other words, we are not interested in reconstructing the interaction function $f$, but only the network structure. At this stage feature ranking can be done with any of the existing feature ranking algorithms, for the purposes of this paper, we consider the two already mentioned algorithms, Random Forest~\cite{BREIM01} and RReliefF~\cite{Kononenko1997}. 

Now, by applying a feature ranking algorithm $R$ to the training data $D_i$ we get feature ranks (importance scores) for node $i$
\begin{equation}
  R(D_i) = (F_{i1}, F_{i2}, \dots , F_{iN}),
\end{equation}
where $F_{ij}$ tells us what is the estimated importance of node $j$ for the node $i$. Note that the values $F_{ij}$ are relative and do not have any physical meaning. By extracting these feature importance scores for all $N$ regression problems, we construct a matrix $F$ of dimension $N \times N$:
\begin{eqnarray}
  \label{fimp}
  F &=& \left(\begin{matrix}
                     F_{11} & F_{12} & \dots  & F_{1N} \\
                     F_{21} & F_{22} &        & \vdots \\
                     \vdots &        & \ddots &  \\
                     F_{N1} & \dots  &        & F_{NN}  
                   \end{matrix}\right) .
\end{eqnarray}
Now, each element $F_{ij}$ in this matrix quantifies how much is the node $i$ important for the node $j$. Our assumption is that higher the value $F_{ij}$ is (relative to other matrix elements) more likely it is that the link $i$---$j$ exists. Hence, we simply extract the reconstructed adjacency matrix $\hat{A}$ from $F$ by setting the threshold $\theta$ and assuming that the links only exist for values of $F$ above the threshold. In general, we can construct $N^2$ different reconstructed adjacency matrices $\hat{A}^{n}$ from $F$ by using each of its elements as a threshold value $\theta_n$:
\begin{equation}
  \label{eq:decision}
  \hat{A}^{n}_{ij}= \begin{cases}
                0; & \text{if $F_{ij}$} \leq \theta_n\\
                1; & \text{if $F_{ij}$} > \theta_n\\
              \end{cases}
\end{equation}
where $n=(1,2,...,N^{2})$.

\paragraph{Measuring reconstruction quality.}  To evaluate how the reconstructed adjacency matrix compares to the real one, we compute a confusion matrix as presented in Table~\ref{table:cm}. 
\begin{table}[h]
  \centering
  \begin{tabular}{l|l|c|c|c}
    \multicolumn{2}{c}{}&\multicolumn{2}{c}{Reconstructed Adj.\ Matrix}&\\
    \cline{3-4}
    \multicolumn{2}{c|}{}&Pred.\ Link&Pred.\ No-Link&\multicolumn{1}{c}{Total}\\
    \cline{2-4}
    \multirow{2}{*}{True Adj.\ Matrix}& Link & TP & FN & TP+FN\\
    \cline{2-4}
    & No-Link & FP & TN & FP+TN\\
    \cline{2-4}
    \multicolumn{1}{c}{} & \multicolumn{1}{c}{Total} & \multicolumn{1}{c}{TP+FP} & \multicolumn{1}{c}{FN+TN} & \multicolumn{1}{c}{$N$}\\
 \end{tabular}
  \caption{Confusion matrix comparing the true adjacency matrix $A$ and the reconstructed one $\hat{A}^{n}$. }
  \label{table:cm}
\end{table}
This confusion matrix tells us much more than the simple accuracy of the reconstruction of all links. For instance, not only we see the number of correctly predicted links and no-links (TP and TN respectively), but also the number of no-links predicted as links (false positives FP) and the number of links predicted as no-links (false negatives FN). We evaluate the performance of the reconstruction in terms of the Sensitivity or True Positive Rate (TPR) and the Fall-Out or False Positive Rate (FPR)~\cite{Fawcett:roc}. These measures are defined as follows:
\begin{equation}
  \mathrm{TPR}=\frac{\mathrm{TP}}{\mathrm{TP}+\mathrm{FN}}, \hspace{1cm} \mathrm{FPR}=\frac{\mathrm{FN}}{\mathrm{FN}+\mathrm{TN}}.
\end{equation}
The TPR tells us the ratio between the correctly predicted links and the total true links, while the FPR is the ratio between the predicted links that are actually no-links and the total number of no-links in the network.

From these two quantities we can further construct the Receiver Operating Characteristic (ROC) curve~\cite{Fawcett:roc} by computing the TPR and the FPR for different thresholds $\theta_{n}$. The ROC curve enable us to evaluate our method without pre-selecting a specific threshold value: It actually includes the results for all possible threshold values. The ROC curve is namely a plot in the TPR and FPR space which presents all network reconstructions that our method produces, and by connecting all the dots in the plot one can compute the Area Under ROC curve (AUC). The larger the AUC (max=1), the better is the method's performance, AUC=1 represents ideal network reconstruction, whereas AUC=0.5 represents reconstruction that is equivalent to random guessing of the link presence. AUC is a measure frequently used in machine learning for binary prediction, and example of which is also our network link prediction.



\section{Results}

In this Section we examine the performance of our reconstruction method. We begin by defining the dynamical system that we will employ (but of course, the method will not use that information). The above formulation of our method is based on discrete-time systems defined via difference equations. However, as we already noted, the method works also for continuous systems provided we can observe and measure the trajectory. Given this, we decided to utilize a discrete-time dynamical system (map) for studying the performance of our method. The key benefit of this is that we need not to worry about measurement resolution and the implications for the precision of derivative estimates. Results presented in this Section are fully applicable to the case of continuous-time dynamical systems as well. In last Section (Discussion) we shall devote more attention to generalization to continuous-time dynamical systems.

For studying the performance of our reconstruction method we choose the logistic map, defined as 
\begin{equation}
  x(t+1) = r x(t) (1 - x(t)) . 
  \label{logistic}
\end{equation}
Logistic map is a textbook example of discrete-time chaotic dynamical system \cite{grassberger1983logi,strogatz2018nonlinear}. The nature of the dynamics, primarily the chaoticity of the behavior, depends on the parameter $r$. Specifically, for $r=4$ the dynamics of Eq.~\ref{logistic} results in chaotic time series behavior for most of the initial conditions. 
To design a dynamical network, we attach a logistic maps to each node $i$ and couple them as done in \cite{masoller2005random,masoller2009synchronizability,Rubido2014a}. The joint equation reads:
\begin{equation}
x_{i}(t+1)=(1-\varepsilon)f(x_{i}(t),r)+\varepsilon\sum\limits^{N}_{j=1;j\neq i}\frac{A_{ij}}{d_{i}}f(x_{j}(t),r),
\label{eq:map}
\end{equation}
where the function $f(x_{i}(t),r\!=\!4)= r x_{i}(t) (1-x_{i}(t))$ stands for the logistic map in the chaotic regime. The parameter $\varepsilon$ denotes the coupling strength between the nodes, and $d_{i}$ the in-degree of each node. We use random networks defined by the link probability ($\rho$) \cite{erdos59,dorogovtsev2008critical}. Each directed link has a probability $\rho$ to be generated. During our study, we will keep the probability to $\rho=0.1$.  

Our method is programmed in Matlab and Python and is available for download at \url{https://github.com/MGrauLeguia/rdn-fr}.
To compute feature rankings we use the Random Forest scikit-learn implementation~\cite{scikit-learn} in Python with $1,000$ trees, we consider square root of all possible features when looking for the best split in a tree, and the RReliefF Matlab implementation (Statistics and Machine Learning Toolbox) with $10$ nearest neighbours. All other parameters of the feature ranking algorithms were the default ones.

\begin{figure*}[htb]
\centering
\includegraphics[width=0.9\linewidth]{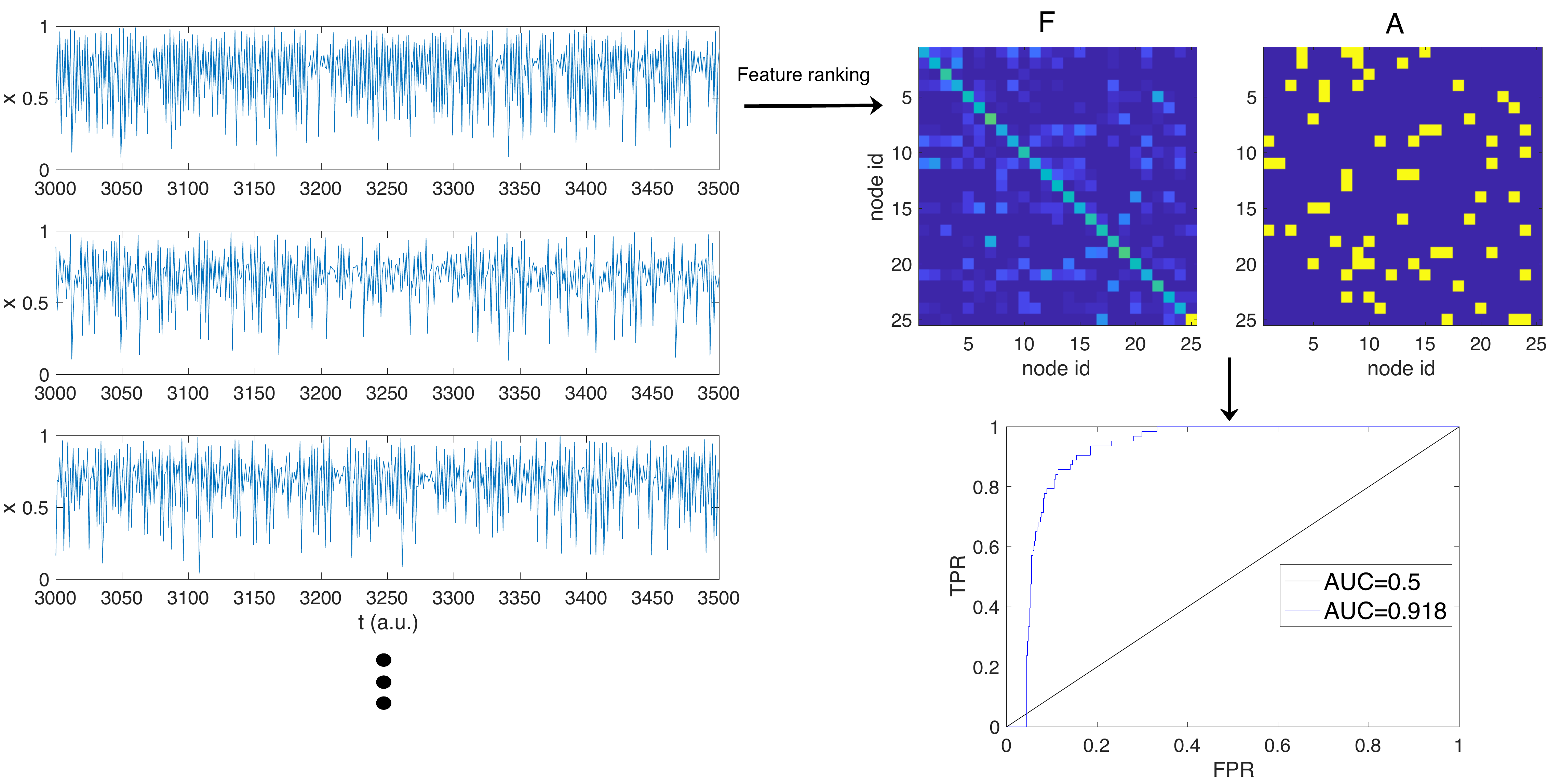}
\caption{A schematic view on how the method works. We first take time series measurements at all the nodes (left-hand side), and use a feature ranking algorithm (e.g., Random Forest or RReliefF) to compute the matrix $F$ containing all the features importance scores. With the selection of a threshold we can get a single solution of the reconstructed network and compare it to the real adjacency matrix $A$ (right-hand side, top). Each threshold value corresponds to single point in the ROC space. However, we evaluate our reconstruction for all possible thresholds by constructing the ROC curve (blue line) and computing the area under it (AUC) (right-hand side, bottom). In this example we used $\varepsilon=0.5$ with a network of $N=25$ nodes, time series with $L=12\,800$ points, and Random Forest for feature ranking. The black line in ROC figure and the area below it denotes the expected behavior for a random reconstruction.}
\label{Fig:example}
\end{figure*}

\subsection{An illustrative example} 

To illustrate our method we examine an example of time series obtained from Eq.~\ref{eq:map}. We consider a random network with $N=25$ nodes and set the coupling strength to $\varepsilon=0.5$. We run it for a random selection of initial conditions and store the obtained time series for each node. The procedure of reconstruction is illustrated in Fig.~\ref{Fig:example}. On the left-hand side the plots show the time series, also to illustrate the nature of signals we are dealing with. On the right-hand side, we show the matrix of the feature importance scores $F$ computed with the Random Forest method and the corresponding true adjacency matrix $A$. We can see that $F$ attains its maximum values along the diagonal corresponding to high self-dependence of all the nodes. This makes sense since at $\varepsilon=0.5$ self-dependence is still high. However, the diagonal elements in $F$ are not taken into account for the calculation of the ROC curve since we do not consider self-loops in the network. Finally, the bottom right part of the figure shows the ROC curve and its corresponding area under it AUC$=0.92$ which we use to measure the performance of the network reconstruction. Procedure is equivalent in case of RReliefF.

\begin{figure*}[htb]
\centering
\includegraphics[width=0.93\linewidth]{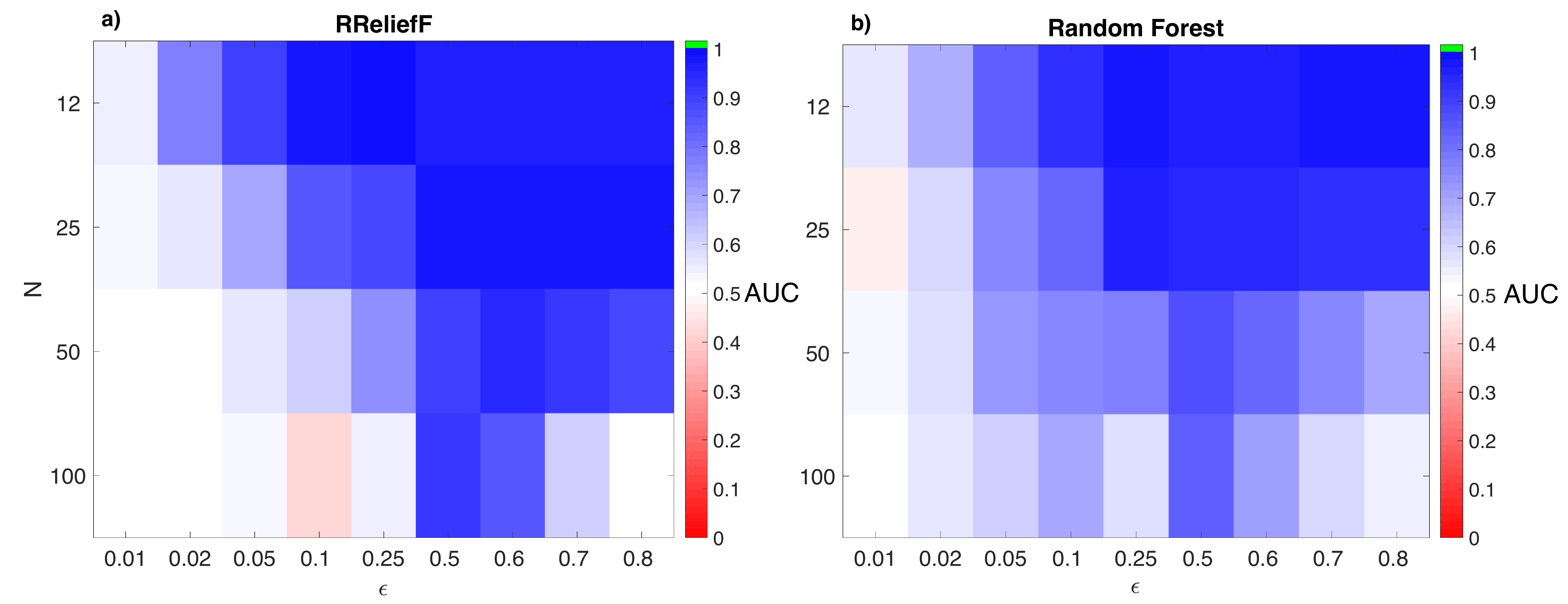}
\caption{Average area under the ROC curve (AUC) across $4$ different realizations of the adjacency matrix $A$ for different network sizes $N$ and coupling strengths $\varepsilon$ using RReliefF (left-hand side) and Random Forest (right-hand side). All input time series comprised $L=12\,600$ data points. With both feature ranking methods we get high performance for high coupling strength and small network size.}
\label{Fig:perf}
\end{figure*}

\begin{figure*}[htb!]
\centering
\includegraphics[width=0.95\linewidth]{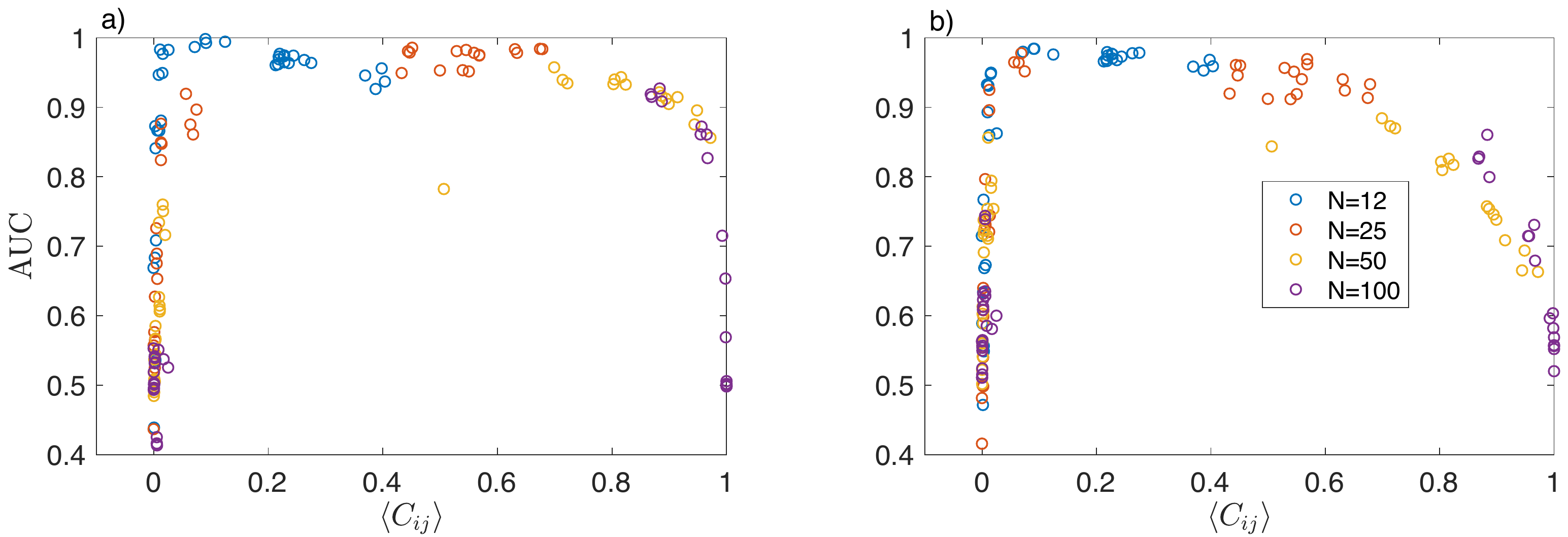}
\caption{Average mean correlation $\langle C_{ij} \rangle$ of the network dynamics plotted against the reconstruction performance (AUC). We used RreliefF ($a$) and Random Forest ($b$) for the feature ranking method using $L=12,800$. We use all combinations of systems size $N$ and coupling strength $\varepsilon$ as in Fig.~\ref{Fig:perf}. Specifically, points with the same system size include all values of $\varepsilon$.}
\label{Fig:coupl.correl}
\end{figure*}

\subsection{Dependence on the size and coupling strength}
\label{res:b}

Next we examine systematically how does the performance depend on the size of the network (number of nodes) and the coupling strength. To this end, we make a grid of parameters ($N$, $\varepsilon$). For each pair (each combination) we draw four different realization of random adjacency matrices and random initial conditions chosen from $[0,1]$. For each of these random realizations, we generate $L=12 800$ data points to which we apply the feature ranking method via both algorithms.

First, we study the performance of the method for a range of coupling strengths $\varepsilon$ and network sizes $N$, while keeping the length of input time series constant and relatively large ($L=12\,800$). In Fig.~\ref{Fig:perf} we present the performance of the method with AUC that is averaged over $4$ independent realizations of the adjacency matrix $A$. We see that at very low coupling strength ($\varepsilon=0.01$) for all network sizes, neither method (RReliefF -- left-hand side, Random Forest -- right-hand side) is able to reconstruct the underlying network---the performance is comparable to the random reconstruction (i.e., AUC $\approx 0.5$).

As the coupling strength increases, Random Forest performs better than RReliefF, especially at large network sizes. This shows that as we increase the network size RReliefF needs higher coupling strength to detect network interactions. This is not the case for Random Forest, with which we find better than random (i.e., AUC$ > 0.5$) reconstruction in areas where RReliefF is performing almost as a random reconstruction. However, around $\varepsilon=0.5$, RReliefF starts to improve. Moreover, the impact of increasing network size on the performance is lower at this coupling strength, and we find very good reconstruction performance for both reconstruction methods at $N=100$ nodes (AUC$_{\textrm{RReliefF}}=0.91$, AUC$_{\textrm{RF}}=0.82$). For $\varepsilon>0.5$ we find very good reconstruction for sizes ranging from $N=12$ to $N=50$, especially for RReliefF, which is still outperforming Random Forest at these high couplings. Finally, at $\varepsilon=0.8$ and $N=100$, neither method is able to correctly detect network interactions and we end up with reconstruction performance close to random.

However, Fig.~\ref{Fig:perf} hides one peculiarity worth examining further. Looking at the case of the largest considered network, $N=100$, we see that for both algorithms performance improves with growing of $\varepsilon$ until approximately $\epsilon=0.5$, but than deteriorates and actually reaches minimum (AUC=0.5) for $\epsilon=0.8$. Intuitive explanation is the following: For small coupling strengths there is not enough interaction among the nodes (logistic maps), so their individual chaoticity prevails, and no useful information can be extracted from the trajectories. In the opposite extreme, for very large coupling strengths the interaction is strong enough to induce very correlated dynamics of nodes/maps (synchronization), and such trajectories also fail to reveal useful information about the underlying network. But, between these two extremes, for intermediate coupling strengths, the interaction might be generating peculiar collective effects that \textit{do reveal} details about the underlying network structure, which are detected by our reconstruction method.

To test this hypothesis we compute the average pair-wise correlation between trajectories $\langle C_{ij} \rangle$. Strongly chaotic trajectories will have (close to) zero $\langle C_{ij} \rangle$, whereas fully regular (synchronized) trajectories will have $\langle C_{ij} \rangle$ equal (or close to) one. In Fig.~\ref{Fig:coupl.correl} we scatter plot the value of $\langle C_{ij} \rangle$ against the value of AUC ($\Omega$) for all the $\varepsilon$'s and realizations. 

And indeed, for small values of $\langle C_{ij} \rangle$, as well as for large values of $\langle C_{ij} \rangle$, the performance is bad for both Random Forest and RReliefF. However, when $\langle C_{ij} \rangle$ is in the intermediate range the performance is good, and it is in fact excellent for a rather wide range of $\langle C_{ij} \rangle$ between roughly 0.1 and 0.9. This indicates that the dynamical regime is intimately related to the ``reconstructibility'' of networks: The reconstruction is clearly the best when the coupling strength is intermediate. There the coupling strength is high enough to reveal important details about its internal structure without making the system fall into a synchronous state. We also note that the system reaches the synchronous state (correlations close to 1) only for the largest system size. This is due to the increase of the average link per node that happens in the larger systems (since we kept the link density constant). This increment of the number links allows the system to be more 'synchronous', and thus exhibit stronger correlations. The performance on most of those cases is random and the non random AUC that we find is due to the initial transients.

\subsection{Dependence on the length of the input time series}

\begin{figure}[htb]
\centering
\includegraphics[width=\linewidth]{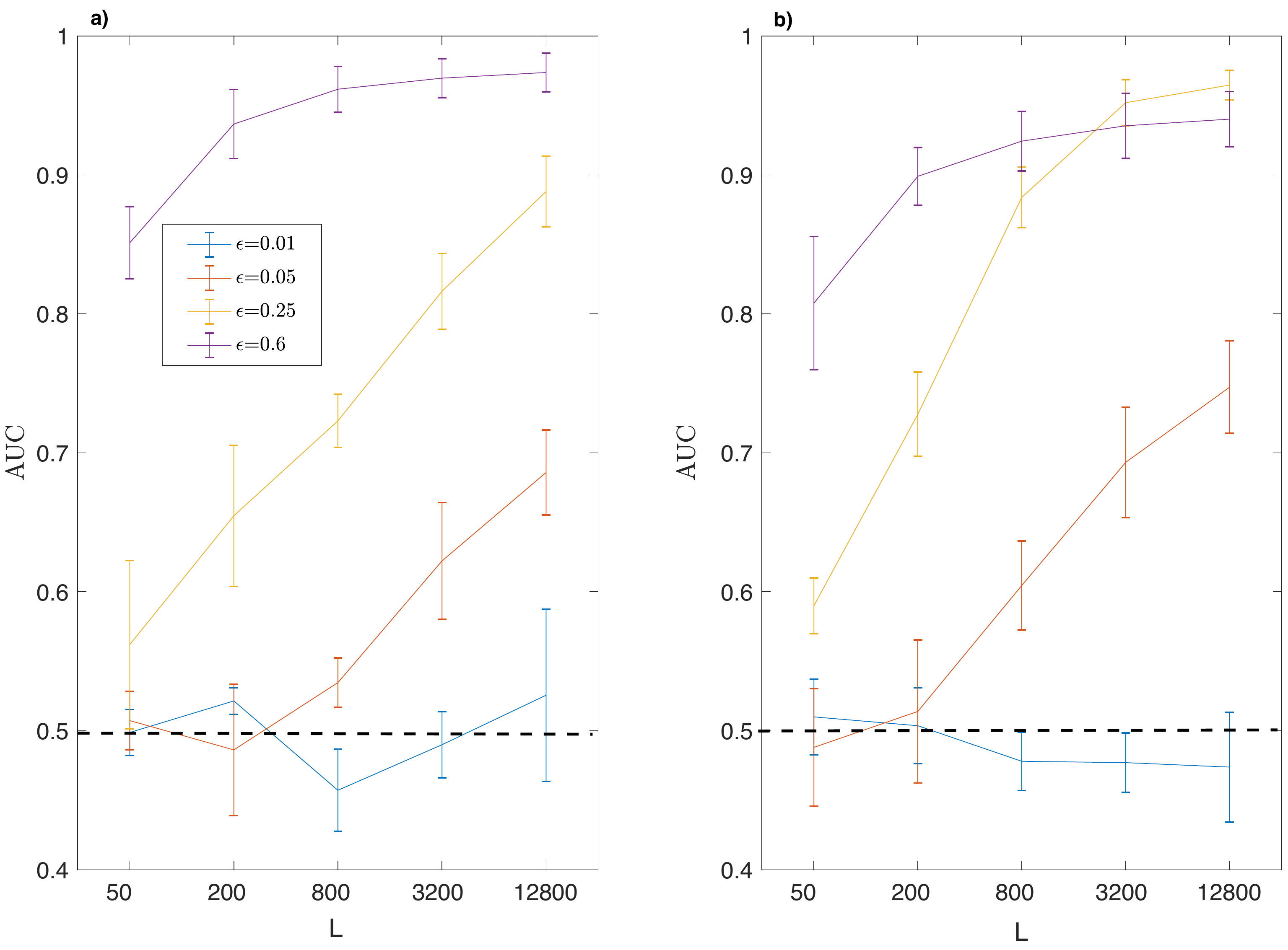}
\caption{Average area under the ROC curve (AUC) as a function of the length of the input time series. We keep the network size fixed to $N=25$, left-hand side presents results for RReliefF and right-hand side for Random Forest. Error bars represent the standard deviation for the four different network realizations. Dashed black line depicts the baseline of the random reconstruction (AUC$ = 0.5$).}
\label{Fig:SizeRR}
\end{figure}

Next we investigate the influence on the input time series length on the performance of the method. In Fig.~\ref{Fig:SizeRR} we present the performance as a function of the time series lengths for both RReliefF (left-hand side) and Random Forest (right-hand side). 
We keep the network size constant $N=25$ and plot the performance for different characteristic coupling strengths $\epsilon$. When the coupling strength is low ($\epsilon=0.01$), increasing the input time series length does not improve the reconstruction performance for any of the two methods and the performance remains close to the one of the random reconstruction. At $\epsilon=0.05$ we start to get better than random at input length $L=800$ for both methods, which perform similarly. Then, at $\epsilon=0.25$, as $L$ increases, the quality of the reconstruction increases at a higher rate for Random Forest than for RReliefF. Finally, at $\epsilon=0.6$, we have a high reconstruction performance even for very short input length $L=50$. Here, adding more input data points does not improve the reconstruction performance dramatically and at $L=800$ the performance only slowly increases.

\begin{figure}[htb]
\centering
\includegraphics[width=\linewidth]{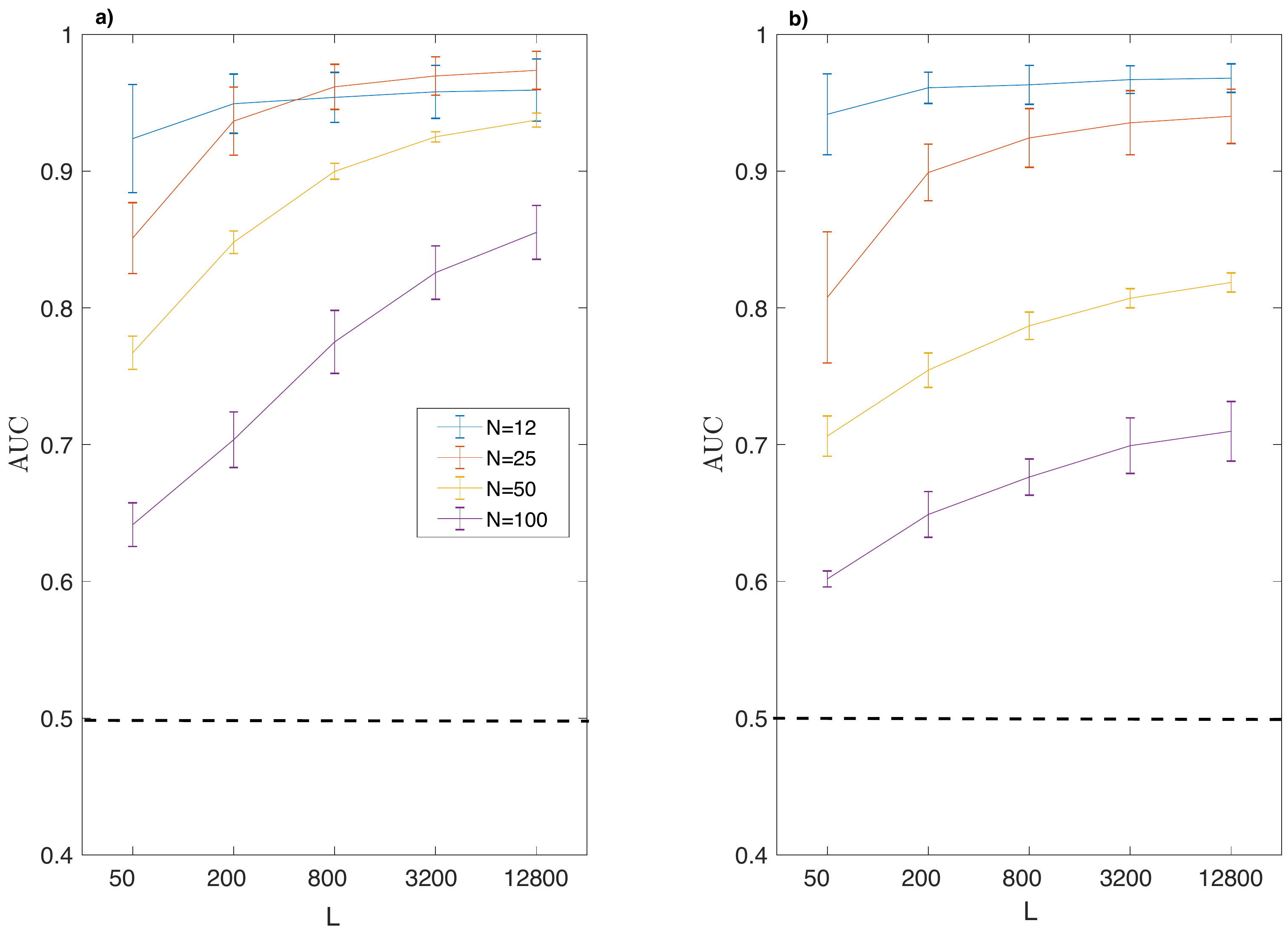}
\caption{Same as in Fig.~\ref{Fig:SizeRR} but keeping the coupling strength constant $\epsilon=0.6$ and plotting results for different network sizes $N$.}
\label{Fig:SizeRF}
\end{figure}

To get a further insight into the reconstruction performance for short input time series lengths, we now keep the coupling strength constant $\epsilon=0.6$ and in Fig.~\ref{Fig:SizeRF} we plot the performance of the reconstruction as a function of $L$ for different network sizes $N$. 
We see that for both feature ranking methods, at $N=12$, the performance is almost perfect and the additional time series data points only slightly improve the reconstruction. We find a similar behavior at $N=25$ where for time series length $L=200$ the performance is already very high (AUC$>0.9$) and only slightly increases with longer time series. As the network size increases ($N=50, 100$), the performance for short input time series $L=50$ decreases significantly. However, our method still performs better than random even for small values of $L$, when $N(N-1)/2 \approx LN$: In other words, when the total number of links, $N(N-1)/2$, to be predicted is similar to the total number of data points, $LN$, we use for the prediction. This suggests that we can use our method even when the available time series is very short.


\subsection{Influence of the noise}

We investigate the impact of noise on the performance of our reconstruction method. We select a realistic set up using observational white noise with zero mean, which, as opposed to the dynamical noise, does not affect the evolution of the network node dynamics. Here, we use the original system from Eq.~\ref{eq:map} and once it is simulated, we add the noise with zero mean and an amplitude of $\sigma$.
\begin{equation}
  \hat{x}_{i}(t) = x_{i}(t) + \xi;
\end{equation}
with $\xi$ as the observational white noise with standard deviation $\sigma$. We now repeat all the computations on $\hat{x}_{i}(t)$ instead of on $x_{i}(t)$. We have to keep in mind that the addition of noise changes both the target value $\hat{x}_{i}(t+1)$ and the features ($\hat{x}_{1}(t)$, $\hat{x}_{2}(t)$,\dots, $\hat{x}_{N}(t)$) are influenced by the observational noise.

\begin{figure}[htb]
\centering
\includegraphics[width=\linewidth]{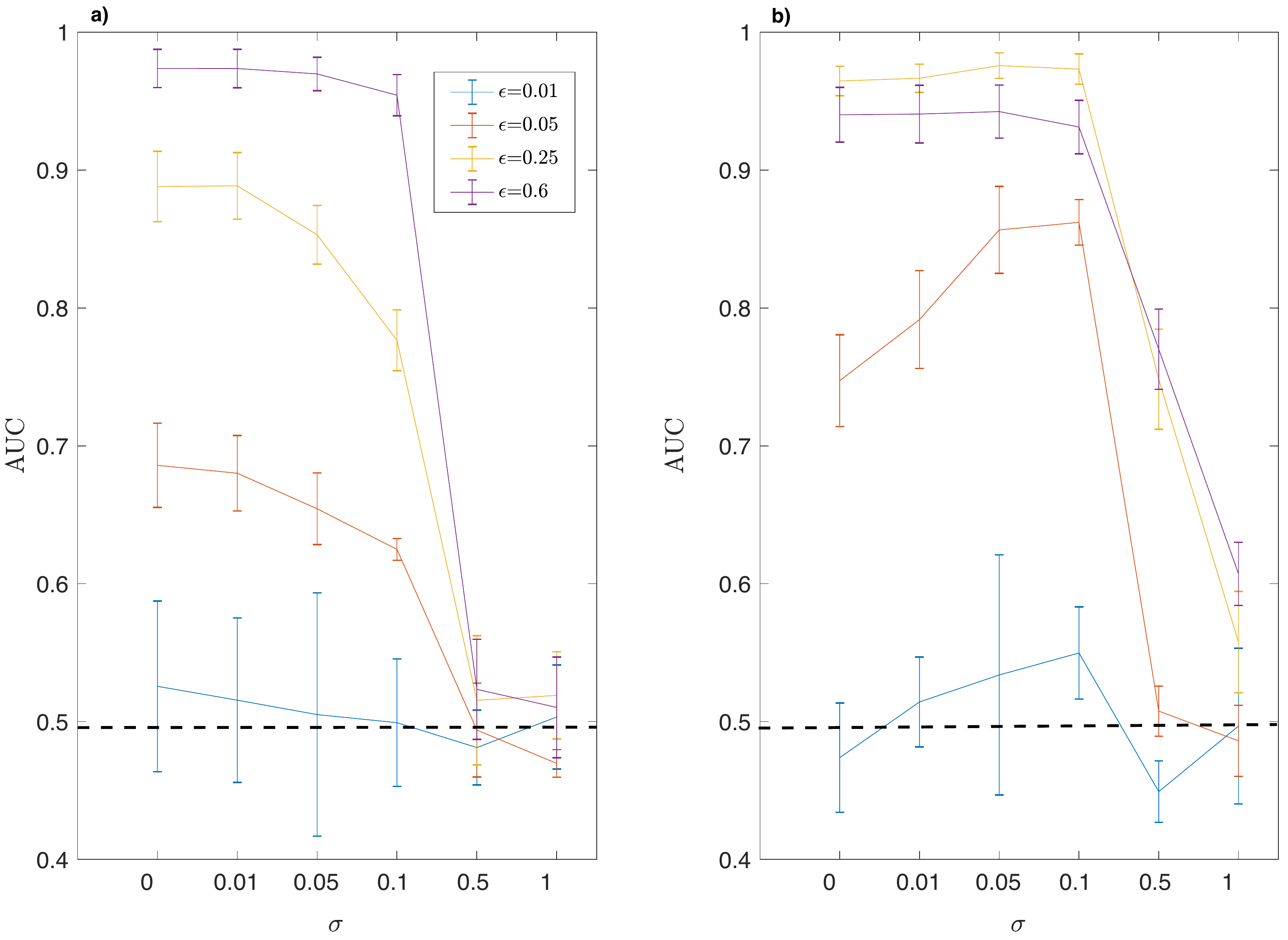}
\caption{Average area under the ROC curve (AUC) as a function of the noise amplitude for different coupling strengths $\epsilon$. We keep the network size constant $N=25$ and present the results for RReliefF on the left-hand side and for Random Forest on the right-hand side. Error bars represent the standard deviation for the four different network realizations. Dashed black line depicts the baseline of the random reconstruction (AUC$ = 0.5$).}
\label{Fig:Noise}
\end{figure}

In Fig.~\ref{Fig:Noise} we present the average area under the ROC curve AUC with RReliefF (left-hand side) and Random Forest (right-hand side) as a function of the amplitude of the Gaussian white noise $\sigma$. We keep $N=25$ and plot the performance for different coupling strengths $\epsilon$. 
At $\epsilon=0.01$, noise does not really influence our results and the reconstruction performance is still close to random. At higher coupling strengths, the change in performance as the noise amplitude grows is lower. For instance, at $\epsilon=0.06$ the performance for both feature ranking methods almost does not change until $\sigma=0.5$. At $\sigma=0.5$, with RReliefF we perform close to random for all the coupling strengths. However, Random Forest is more robust as even at $\sigma=0.5$ we have a better than random reconstruction for $\epsilon=0.25, 0.6$. This effect can be explained due to the fact that noise is added to both sides of Eq.~\ref{eq:map}. Finally, at $\sigma=1$ the performance decreases until we get close to random performance. We have to keep in mind that the two last noise amplitudes considered ($\sigma=0.5,1$) have a similar or higher amplitude to the original time series amplitude, which is an overestimation of real world examples of noise levels. Therefore, for both feature ranking methods (and especially for Random Forest, which is itself known to be robust to noise) our method is very robust to observational noise.

\subsection{Application to more complex dynamical systems}
While the logistic map served us as a simple model to test the general performance of our method, in this section we show that our method is also useful for more complicated dynamical systems. Specifically, for the function $f$ in Eq.~\ref{eq:map} we now consider the Ikeda complex map map~\cite{ikeda}: 
\begin{eqnarray}
  x_{i}(j+1)&=&1+u(x_{i}(j)\cos t_{i}(j)-y_{i}(j)\sin t_{i}(j))\\
  y_{i}(j+1)&=& u(x_{i}(j) \sin t_{i}(j)+y_{i}(j)\cos t_{i}(j))
\end{eqnarray}
with $t_{i}(j)=0.4-\frac{6}{1+x^{2}_{i}(j)+y^{2}_{i}(j)}$. The parameter $u$ controls the dynamics and for $u>0.6$ a chaotic attractor is generated.
\begin{figure}[htb]
\centering
\includegraphics[width=\linewidth]{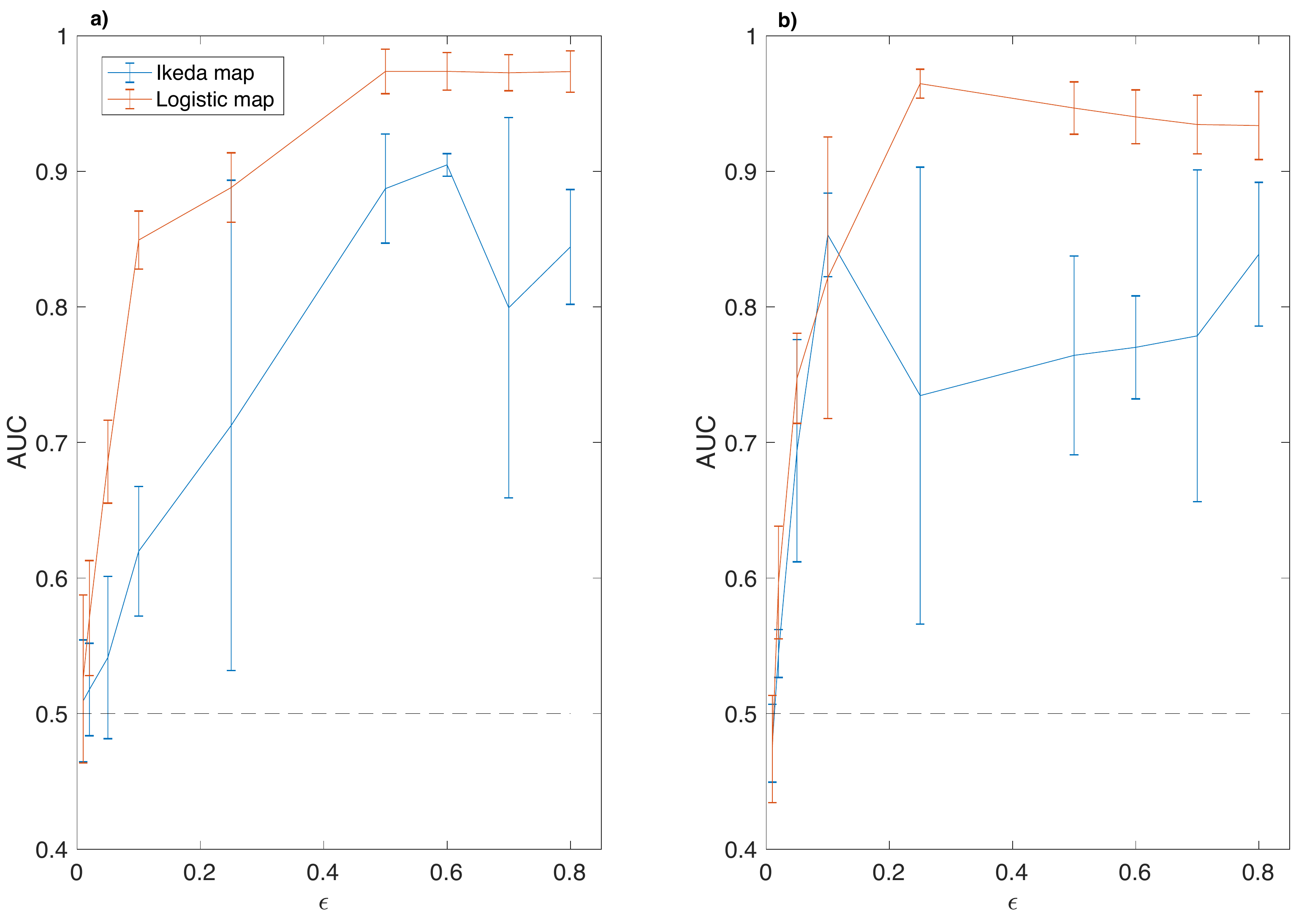}
\caption{Average area under the ROC curve (AUC) as a function of the coupling strengths $\varepsilon$ for the RReliefF on the left-hand side and for Random Forest on the right-hand side. We keep the network size constant with $N=25$ and use $L=12\,800$. The Ikeda map was computed using $u=0.9$ whereas the logistic map was computed using $r=4$. We used the same values of $\varepsilon$ as in Fig.~\ref{Fig:perf}. Error bars represent the standard deviation for the four different network realizations. Dashed black line depicts the baseline of the random reconstruction (AUC=0.5).}
\label{Fig:Ikeda}
\end{figure}
In the following we will keep $u=0.9$ to ensure we are in a region with sufficient complex dynamics. We now simulate the network of coupled Ikeda maps and in each node get a time series with complex values. We select only the real part of the time series, and use this data as an input to our network reconstruction method. Therefore, in this section, we are dealing with more complex dynamics and we also no longer have complete observability of our system (the imaginary part is not known to our method).

In Fig.~\ref{Fig:Ikeda}, we present the averaged performance of our method using coupled Ikeda maps (blue) and Logistic maps (red) for $4$ different realizations of the adjacency matrix $A$. We kept the network size to $N=25$ and the length of the time series to $L=12\,800$. As in Fig.~\ref{Fig:perf}, in Fig.~\ref{Fig:Ikeda}(a) we used RReliefF and at Fig.~\ref{Fig:Ikeda}(b) we used Random Forest.

Similarly to the results of the coupled logistic maps, for lower coupling strengths (RReliefF:~$\varepsilon \leq 0.02$, Random Forest:~$\varepsilon \leq 0.05$) the reconstruction with Ikeda maps is marginally better than random (AUC $\leq 0.5$). However, for higher coupling strengths the performance improves and we get good reconstructions (AUC $\approx 0.8$ or more). Nevertheless, the performance is lower than the one we have with logistic maps (AUC $\approx 0.95$). This is expected as we are comparing a complex map with a one dimensional map, and only partial information (the real part of the Ikeda map) of the complex map is used.

Interestingly, our methods with RReliefF and Random Forest perform best for the Ikeda maps at different coupling strengths. Specifically, with Random Forest it starts yielding good reconstructions even for very small coupling strengths, meaning it can extract useful information even from very weak self-organization. In contrast, with RReliefF it needs a considerable coupling strength for achieving higher AUC, but it then outperforms the Random Forest version. Moreover, we observe that the AUC scores vary quite considerably at each realization. This strong variability is not observed when we use the logistic maps and could be caused by the fact that we now only have partial observability and we are only using partial information of the system (we compute the scores using only the real part of the Ikeda maps).

\subsection{Comparison with Partial Mutual Information from Mixed Embedding (PMIME)}
Finally, we check how our method compares to existing data-driven methods. For comparison we selected the Partial Mutual Information from Mixed Embedding (PMIME) method~\cite{Koutlis2016a} that is built around an entropy-based measure that can detect directionality of links~\cite{Kugiumtzis2013}. The method, like ours, does not require any strong assumptions on the nature of the reconstructed system, and its usefulness has been demonstrated on Mackey-Glass delay differential equations and neural mass models. PMIME requires setting a termination criterion to avoid false positives, but in order to compute the ROC statistics and directly compare PMIME with our method, we set this parameter to zero and evaluate the resulting connectivity matrix in the same manner as the connectivity matrix produced by our method. All the other PMIME parameters were set to their default values.

\begin{figure}[htb]
\centering
\includegraphics[width=\linewidth]{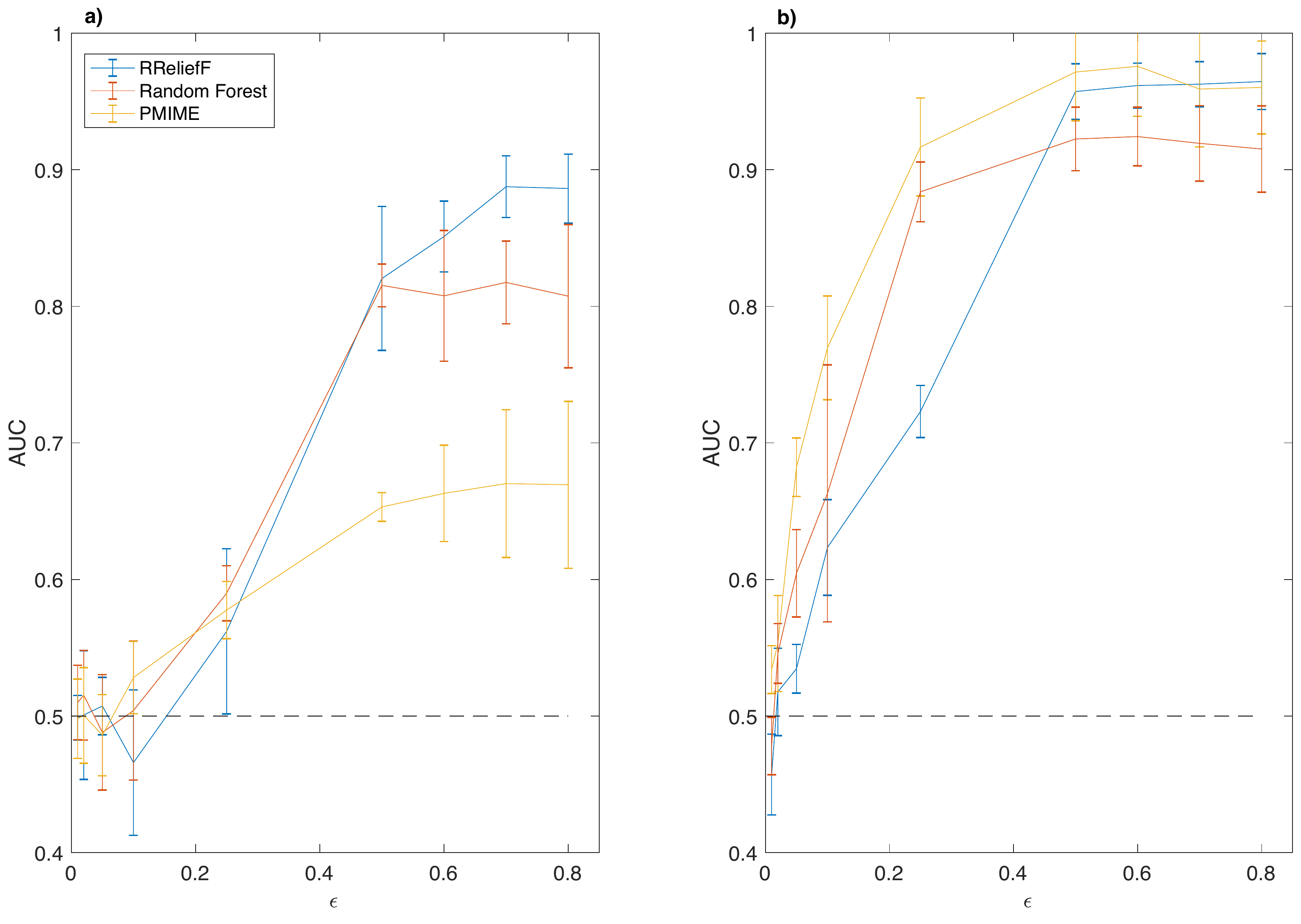}
\caption{Average area under the ROC curve (AUC) as a function of the coupling strengths $\varepsilon$ for the RReliefF, Random Forest and PMIME using $L=50$ (a) and $L=800$ (b). We keep the network size constant with $N=25$. Error bars represent the standard deviation for the four different network realizations. Dashed black line depicts the baseline of the random reconstruction (AUC=0.5).}
\label{Fig:PMIME}
\end{figure}

In Fig.~\ref{Fig:PMIME} we show the performance PMIME and both variants of our method as a function of the coupling strength for two different time series lengths $L=50$ (Fig.~\ref{Fig:PMIME}a) and $L=800$ (Fig.~\ref{Fig:PMIME}b) keeping the network size at $N=25$. At longer time series ($L=800$) and low and intermidiate coupling strengths ($\varepsilon < 0.5$) PMIME is performing better, especially in comparison to the RReliefF variant. For $\varepsilon \geq 0.5$, both variants of our method and PMIME are performing similarly with an almost perfect AUC. For shorter time series ($L=50$), both variants of our method are outperforming PMIME when the coupling strength is big enough ($\varepsilon > 0.25$). It is known that entropy-based methods need long time series for producing reliable estimates, therefore the decrease of performance is expected. On the other hand, the decrease of performance of both variants of our method is much smaller and the RReliefF variant scores an AUC $\approx0.9$ with only $L=50$ points.

\begin{figure}[htb]
\centering
\includegraphics[width=\linewidth]{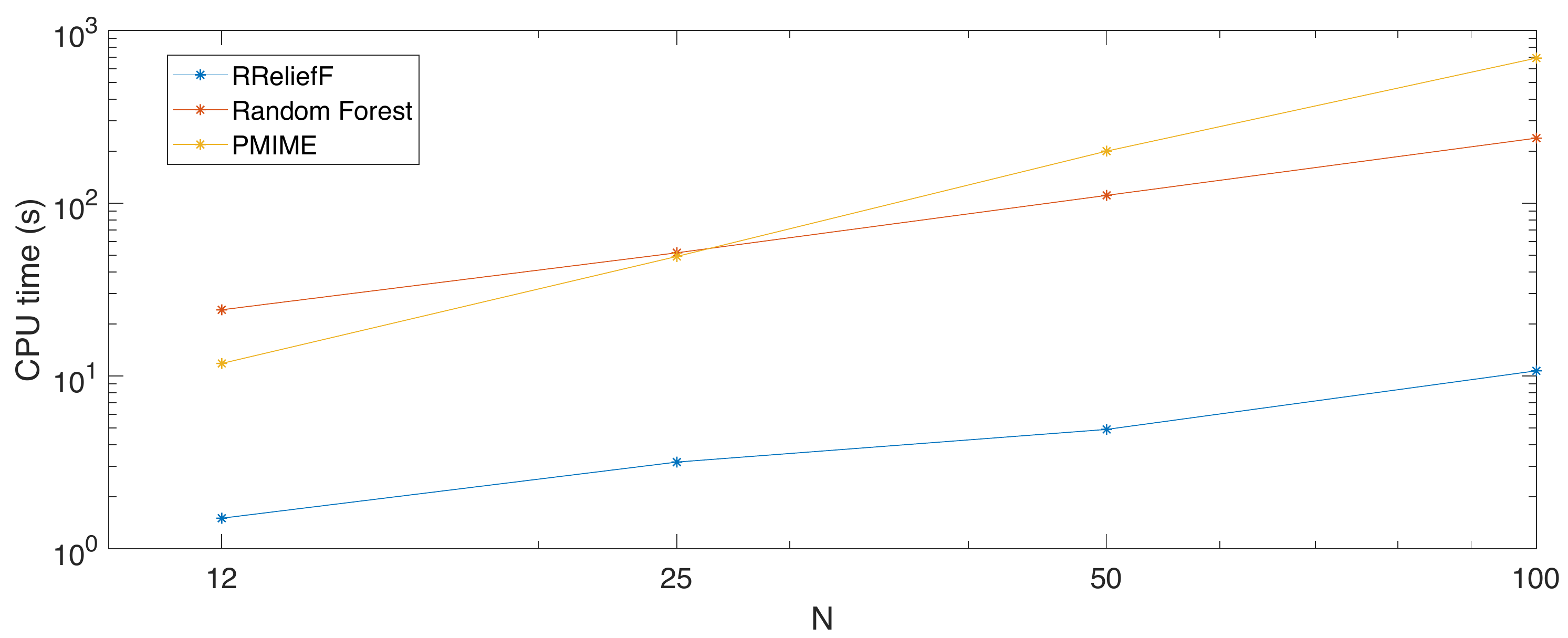}
\caption{CPU times (log scale) of both variants of our method (RReliefF and Random Forest) and PMIME as a function of the system size using time series of length $L=200$.}
\label{Fig:Time}
\end{figure}

Furthermore, we investigated the scaling properties of PMIME and our methods. In Fig.~\ref{Fig:Time} we show CPU times in log scale of all three algorithms as a function of the system size $N$ using time series of length $L=200$. We observe that both RReliefF and Random Forest scale approximately with the same slope of $m\approx1$ suggesting an approximately linear scaling factor with $N$. For PMIME, on the other hand, the slope suggests an approximately quadratic scaling factor with $N$. The plot also shows that at these system sizes the RReliefF variant is by an order of magnitude faster than the Random Forest variant. The linear scaling with the size of the reconstructed system of our method makes them suitable for inferring larger networks in practice.

\section{Discussion}

We designed a novel method for reconstructing (inferring) networks of dynamical units from observation (measurements) of their node trajectories. It is based on feature ranking, a common methodology of machine learning. By ranking the ``features", which are the values of the trajectories of other nodes, we can extract information on what other nodes are most likely to be connected with the considered node. We test the performance of our method using networks of coupled logistic maps, and obtain good results for a range of coupling strengths and network sizes. Also, our method is able to perform well even for relatively short trajectories and it's fairly robust to noise.

The key property of our method is that it requires no assumption on the knowledge of interaction functions or the dynamical model of the network, and that it makes no hypotheses on the nature of the available trajectories (data). We consider this to be an important aspect when selecting a reconstruction method for a practical application, since most, though not all, of the similar methods in the current (physics) literature make assumptions about above mentioned details that can sometimes be rather strong. So, while our method is not based on physical insight into the collective dynamics, it is immediately applicable to practically any complex dynamical system (of physical interest or otherwise), requiring no prior knowledge about system's internal details. 

We demonstrate the wide applicability of our method also with experiments on more complex dynamical systems with only partial observability---on the Ikeda maps. While the performance of the method is, as expected, lower than on the simpler logistic maps, we can still get good reconstructions for higher coupling strengths. 

Given a large number of available network reconstruction methods a question of their comparison arises. However, comparing them is not trivial, since various methods depart from different hypotheses and knowledge about the system, which makes their merits harder to compare. Another distinction is also between what different methods are reconstructing: interaction functions and network structure, or just the structure? Hence, our method can be meaningfully compared only to methods that: (i) reconstruct only the structure, (ii) make no assumptions on interaction functions, and (iii) rely on discrete measurements of dynamical trajectories of nodes. We therefore compare our method to PMIME, which, like our method, also reconstructs the structure and does not require any knowledge of the system. The results show that with PMIME we can obtain comparable or sometimes even better reconstructions for longer time series. However, for shorter time series, which we frequently come across in practice, our method outperforms PMIME. An even more important advantage of our method turns out to be its scalability; while CPU times of our method grow roughly linearly with the size of the system, PMIME grows roughly quadratically. Our method is therefore in practice applicable to much larger systems than PMIME.

Still, our method does have some limitations. Specifically, the performance seems to deteriorate as the system size increases. For reconstructing large networks, long observation trajectories are needed for accurate reconstruction. This could hinder the applicability to large systems that abound in applications. But despite this, for certain dynamical regimes (ranges of coupling strength), the performance remains good independently of system size. This represents a hope for applications to large systems. In contrast, for other dynamical regimes (too small or too strong coupling), the performance worsens. This cannot be helped, since in those dynamical regimens (full chaos or full regularity) the system reveals nothing about its underlying structure. Also, while our method in general reacts well to the noise, too excessive noise deteriorates the performance. On the other hand, as discussed above, our method is not limited by a any assumptions about prior knowledge on the system under scrutiny. 

Another issue revolves around the generalization to continuous systems. This paper's presentation was based on maps, but with minor modifications, the method can be generalized to continuous-time dynamical systems (defined via differential rather than difference equations). To this end, one would need to replace the right-hand side of Eq.~\ref{eq:featureranking} with $\dot{x_i}(t)$, the time derivative (rate of change) of the state of a node $i$ at time $t$. Since we assume total observability of the system, one can compute the derivatives of the state trajectories numerically and use the same regression and feature ranking algorithms to infer the network structure. However, further empirical evaluation of the proposed method is needed in such a setting. While the task of numerical differentiation is known to be unstable, recent methods based on reformulating the task of numerical differentiation as an optimization problem and using regularization methods have been proposed to obtain accurate estimates of the derivatives from noisy and sparse trajectories~\cite{ramm2001,hanke2001}.

The core purpose of this paper is to present the concept and theoretical underpinnings for our method. Next step is to apply it to real systems in two steps. First, one should use a scenario where the ground truth is available to test the performance in this setting. Second, our method can be used to actually reveal the structure of real networks whose structure is still unknown. We envisage one such possible application in brain network inference from electroencephalographic recordings, where the interest is the functional brain connectivity obtained without knowing (or inferring) the interaction function(s). Also, data on gene expression can be used in a similar way to reconstruct gene regulation networks under given circumstances. Our method is likely to serve as a valuable alternative here, at least to some extent, but quantifying how much will require additional work. We here also mention that our method requires no interference with the systems (such as e.g.\ random resets), and is as such non-invasive. However, such interference could further improve its performance.

We close the paper with discussing the possible avenues of future work. That primarily includes improving this method by combining the rankings obtained with different feature ranking methods. Here, we have only compared the performance of our methods when using Random Forest \emph{or} RReliefF algorithms for ranking. In further work, one can also combine an ensemble of rankings obtained with different feature ranking algorithms~\cite{prati2012}. Furthermore, we focus here solely on ranking of features. Instead, one could explore the use of the predictive models that are learned using Random Forest or other supervised learning algorithms. In particular, the predictive models can reveal various relevant aspects of the mechanisms of interactions between network nodes, in many real-world cases when these mechanisms are not known.

\section*{Acknowledgements}
The study was financially supported by the Slovenian Research Agency through research core funding of programs P1-0383, Complex Networks (Z.L.), P2-0103 (B.Ž.) and P5-0093 (L.T.), as well as projects J5-8236 (Z.L.) and N2-0056, Machine Learning for Systems Sciences (L.T.). M.G.L.\ and Z.L.\ acknowledge founding from the EU via H2020 Marie Sklodowska-Curie project COSMOS, (Grant No.\ 642563), and B.Ž.\ acknowledges founding from the EU via H2020 projects SAAM (Grant No.\ 769661) and RESILOC (Grant No.\ 833671). 

\bibliography{iopart-num}
\bibliographystyle{plain}
\end{document}